\documentclass[11pt,english]{article}%
\usepackage[T1]{fontenc}
\usepackage[latin9]{inputenc}
\usepackage{color}
\usepackage{amsmath}
\usepackage{graphicx}
\usepackage{amssymb}
\usepackage{amsfonts}
\usepackage{a4wide}
\usepackage{babel}%
\providecommand{\U}[1]{\protect\rule{.1in}{.1in}}
\newtheorem{theorem}{Theorem}

\newtheorem{example}[theorem]{Example}

\begin{document}

\title{On three duality results}
\author{M. D. Voisei\thanks{Towson University, U.S.A., email:
\texttt{mvoisei@towson.edu}.}~ and C. Z\u{a}linescu\thanks{University
{}``Al.I.Cuza'' Ia\c{s}i, Faculty of Mathematics, Romania, email:
\texttt{zalinesc@uaic.ro}.}}
\date{}
\maketitle

The aim of this short note is to give counterexamples to two results by D. Y.
Gao \cite[Th.\ 16]{Gao/Sherali:07}, \cite[Th.\ 2]{Gao:07} and to improve a
related result by S.-C. Fang, D. Y. Gao, R.-L. Sheu and S.-Y. Wu
\cite[Th.\ 3]{FGSW}.

\section{Counterexamples to \cite[Th.\ 16]{Gao/Sherali:07}, \cite[Th.\ 2]%
{Gao:07}}

On \cite[page 298]{Gao/Sherali:07} the authors consider the problem

{}``$\min\left\{  P(x)=\frac{1}{2}x^{T}Ax-f^{T}x\ :\ \tfrac{1}{2}x^{T}%
Cx\leq\lambda,\ x\in\mathbb{R}^{n}\right\}  $.$\quad$ (8.156)''

\noindent{}``...where $A$ and $C$ are two symmetrical matrices in
$\mathbb{R}^{n\times n}$, $f\in\mathbb{R}^{n}$ is a given vector, and
$\lambda\in\mathbb{R}$ is a given constant'', and continue on the following
page with: {}``On the dual feasible space

$\mathcal{V}_{k}^{\ast}=\{\varsigma\in\mathbb{R}\mid\varsigma\geq0,$
$\det(A+\varsigma C)\neq0\}$

\noindent and the canonical dual problem (8.155) can be formulated as (see {[}50]):

$\max\left\{  P^{d}(\varsigma)=-\tfrac{1}{2}f^{T}(A+\varsigma C)^{-1}%
f-\lambda\varsigma\ :\ \varsigma\in\mathcal{V}_{k}^{\ast}\right\}  .\quad$ (8.158)''

{}``The following result was obtained recently.

\bigskip{}

\textbf{Theorem 16 (Gao {[}50])} Suppose that the matrix $C$ is positive
definite, and $\overline{\varsigma}\in\mathcal{V}_{a}^{\ast}$ is a critical
point of $P^{d}(\varsigma)$. If $A+\overline{\varsigma}C$ is positive
definite, the vector

$\overline{x}=(A+\overline{\varsigma}C)^{-1}f$

\noindent is a global minimizer of the primal problem (8.156). However, if
$A+\overline{\varsigma}C$ is negative definite, the vector $\overline
{x}=(A+\overline{\varsigma}C)^{-1}f$ is a local minimizer of the primal
problem (8.156).''

\bigskip{}

In the previous statement $\mathcal{V}_{a}^{\ast}=[0,+\infty)$ (see
\cite[p.\ 297]{Gao/Sherali:07}) while reference {[}50] is our reference
\cite{Gao:05}. This first result we are interested in is cited in
\cite{Gao/Sherali:07} as being published in \cite{Gao:05}; however, we could
not find its statement in \cite{Gao:05}. The following is a counterexample for
\cite[Th.\ 16]{Gao/Sherali:07}.

\begin{example}
\label{ex-th16-gs}Consider%
\[
A=\left[
\begin{array}
[c]{cc}%
-2 & -1\\
-1 & -3
\end{array}
\right]  ,\quad C=\left[
\begin{array}
[c]{cc}%
1 & 0\\
0 & 1
\end{array}
\right]  ,\quad f=\left[
\begin{array}
[c]{c}%
-1\\
-1
\end{array}
\right]  ,\quad\lambda=\frac{1}{2}.
\]
Then $P^{d}(y)=-\frac{1}{2}y-\frac{1}{2}\frac{2y-3}{y^{2}-5y+5}$ and
$(P^{d})^{\prime}(y)=-\frac{1}{2}\frac{\left(  y-2\right)  ^{2}}{\left(
y^{2}-5y+5\right)  ^{2}}\left(  y-1)(y-5\right)  $. Hence the set of critical
points of $P^{d}$ is $\{1,2,5\}$ all contained in $\mathcal{V}_{a}^{\ast}$.
For $\overline{y}=1$ we have that $A+\overline{y}C=\left(
\begin{array}
[c]{cc}%
-1 & -1\\
-1 & -2
\end{array}
\right)  $ is negative definite and $\overline{x}=\left(  A+\overline
{y}C\right)  ^{-1}f=\left[
\begin{array}
[c]{cc}%
1 & 0
\end{array}
\right]  ^{T}$.

Let $\mathcal{U}_{0}:=\{(\cos t,\sin t)^{T}\mid t\in(-\pi,\pi)\}$ be a subset
of the admissible set $\mathcal{U}=\{x\in\mathbb{R}^{2}\mid\left\Vert
x\right\Vert \leq1\}$ and
\[
f(t):=P((\cos t,\sin t)^{T})=-1-\cos t\sin t-\tfrac{1}{2}\sin^{2}t+\cos t+\sin
t=-(3+\cos t-2\sin t)\sin^{2}\tfrac{1}{2}t,
\]
$t\in\mathbb{R};$ hence $P(\overline{x})=f(0)=0=P^{d}(\overline{y})$.
According to the previous theorem $\overline{x}$ should be local minimizer of
$P$ on $\mathcal{U}$, in contradiction to the fact that $t=0$ is a strict
local maximum point of $f$ (see Figure 1).
\end{example}

\begin{center}

\includegraphics[bb=0 0 300 200,
,scale=0.5]{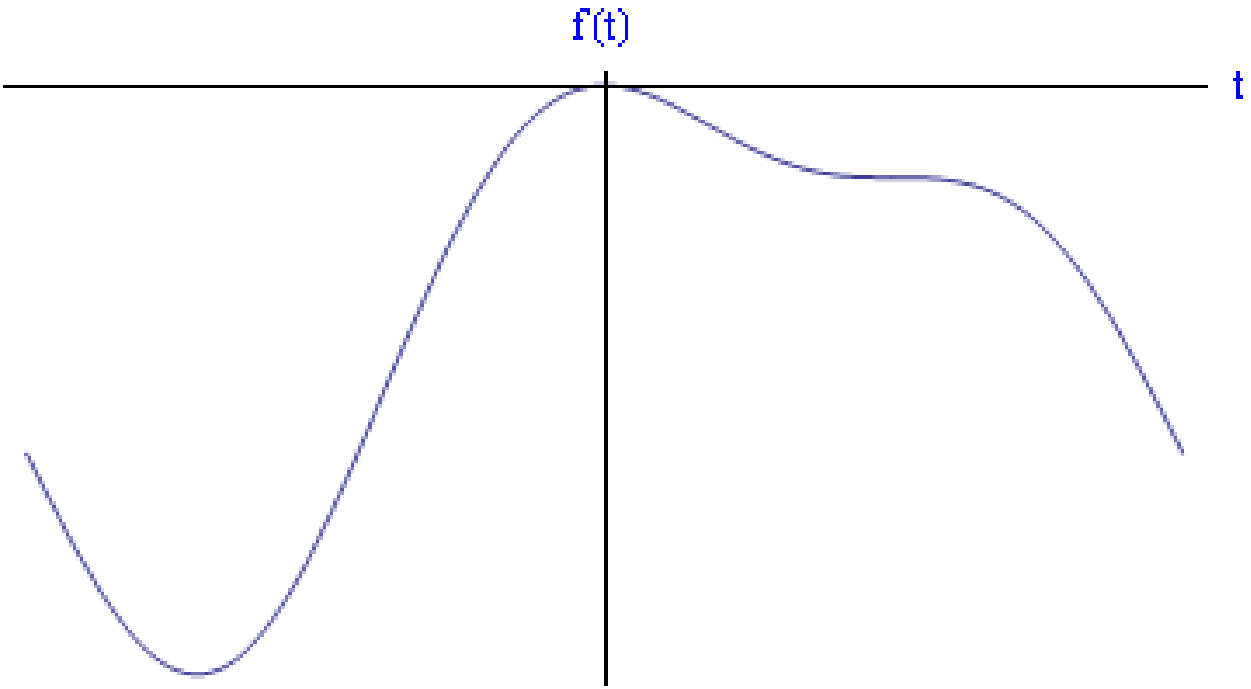}

Figure 1.
\end{center}

\strut

Our attention turns to the problem considered in \cite{Gao:07}

{}``$(\mathcal{P})\ :\ \ \min\left\{  P(x)=U(\Lambda(x))+Q(x)\ :\ x\in
\mathcal{R}^{n}\right\}  \quad$ (5)''

\noindent where {}``$Q(x)=\tfrac{1}{2}x^{T}Ax-c^{T}x$ is a quadratic function,
$A=A^{T}\in\mathcal{R}^{n\times n}$ is a given symmetric matrix'',
$c\in\mathcal{R}^{n}$, and the so called {}``geometrical operator
$\Lambda:\mathcal{R}^{n}\rightarrow\mathcal{R}^{1+n}$ and the associated
canonical function $U$ can be introduced as following:

$y=\Lambda(x)=\left(
\begin{array}
[c]{c}%
\xi(x)\\
\epsilon(x)
\end{array}
\right)  =\left(
\begin{array}
[c]{c}%
\tfrac{1}{2}\left\vert Bx\right\vert ^{2}-\alpha\\
\left\{  x_{i}^{2}-\ell_{i}\right\}
\end{array}
\right)  \in\mathcal{R}^{1+n},$

$U(y)=\tfrac{1}{2}\xi^{2}+\Psi(\epsilon)\quad$ (3)

\noindent where

$\Psi(\epsilon)=\left\{
\begin{array}
[c]{cl}%
0 & \text{if }\epsilon\leq0,\\
+\infty & \text{otherwise.}%
\end{array}
\right.  \quad$ (4)''

\noindent Here {}\textquotedblleft$B\in\mathcal{R}^{m\times n}$ is a given
matrix and $\alpha>0$ is a given parameter\textquotedblright\ while
$\ell=\{\ell_{i}\}\in\mathcal{R}^{n}$, $\ell_{i}\geq0$. {}\textquotedblleft
The notation $\left\vert x\right\vert $ used in this paper denotes the
Euclidean norm of $x$\textquotedblright.

{}``The canonical dual problem of $(\mathcal{P})$ can be proposed as the following

$(\mathcal{P}^{d})$ : $\operatorname*{sta}\left\{  P^{d}(\varsigma
,\sigma)=-\tfrac{1}{2}c^{T}\left[  G(\varsigma,\sigma)\right]  ^{-1}%
c-\tfrac{1}{2}\varsigma^{2}-\alpha\varsigma-\ell^{T}\sigma\ :\ (\varsigma
,\sigma)^{T}\in\mathcal{S}_{a}\right\}  .\quad$ (11)

Here {}``$G(\varsigma,\sigma)$ is a symmetrical matrix, defined by

$G(\varsigma,\sigma)=A+\varsigma B^{T}B+2\operatorname*{Diag}(\sigma
)\in\mathcal{R}^{n\times n}$, $\quad$(9)

\noindent and $\operatorname*{Diag}(\sigma)\in\mathcal{R}^{n\times n}$ denotes
a diagonal matrix with $\{\sigma_{i}\}$ $(i=1,2,\ldots,n)$ as its diagonal
entries'' while {}``$\mathcal{S}_{a}=\big\{  y^{\ast}=\big( \textstyle%
\begin{array}
[c]{c}%
\varsigma\\
\sigma
\end{array}
\big)  \in\mathcal{R}^{1+n}\mid\varsigma\geq-\alpha,\ \sigma\geq0,\ \det
G(\varsigma,\sigma)\neq0\big\}  .\quad$ (10)''

One continues with {}``we need to introduce some useful feasible spaces:

$\mathcal{S}_{a}^{+}=\left\{  (\varsigma,\sigma)^{T}\in\mathcal{S}_{a}\mid
G(\varsigma,\sigma)\text{ is positive definite}\right\}  $,$\quad$ (16)

$\mathcal{S}_{a}^{-}=\left\{  (\varsigma,\sigma)^{T}\in\mathcal{S}_{a}\mid
G(\varsigma,\sigma)\text{ is negative definite}\right\}  $.$\quad$ (17)

\bigskip{}

\textbf{Theorem 2} (Triality Theorem). Suppose that the vector $\overline
{y}^{\ast}=(\overline{\varsigma},\overline{\sigma})^{T}$ is a KKT point of the
canonical dual function $P^{d}(y^{\ast})$ and $\overline{x}=[G(\overline
{\varsigma},\overline{\sigma})]^{-1}c$.

If $\overline{y}^{\ast}=(\overline{\varsigma},\overline{\sigma})^{T}%
\in\mathcal{S}_{a}^{+}$, then $\overline{y}^{\ast}$ is a global maximizer of
$P^{d}$ on $\mathcal{S}_{a}^{+}$, the vector $\overline{x}$ is a global
minimizer of $P$ on $\mathcal{X}_{a}$, and

$P(\overline{x})=\min\limits_{x\in\mathcal{X}_{a}}P(x)=\max\limits_{y^{\ast
}\in\mathcal{S}_{a}^{+}}P^{d}(y^{\ast})=P^{d}(\overline{y}^{\ast}).\quad$ (18)

If $\overline{y}^{\ast}\in\mathcal{S}_{a}^{-}$, on the neighborhood
$\mathcal{X}_{o}\times\mathcal{S}_{o}\subset\mathcal{X}_{a}\times
\mathcal{S}_{a}$ of $(\overline{x},\overline{y}^{\ast})$, we have that either

$P(\overline{x})=\min\limits_{x\in\mathcal{X}_{o}}P(x)=\min\limits_{y^{\ast
}\in\mathcal{S}_{o}}P^{d}(y^{\ast})=P^{d}(\overline{y}^{\ast})\quad$ (19)

\noindent holds, or

$P(\overline{x})=\max\limits_{x\in\mathcal{X}_{o}}P(x)=\max\limits_{y^{\ast
}\in\mathcal{S}_{0}}P^{d}(y^{\ast})=P^{d}(\overline{y}^{\ast}).\quad$ (20)''

\bigskip{}

Recall that {}``$\mathcal{X}_{a}=\{x\in\mathcal{R}^{n}\mid\ell^{l}\leq
x\leq\ell^{u}\}$ is a feasible space'' and {}``we assume without loss of
generality that $\ell^{u}=-\ell^{l}=\ell^{\frac{1}{2}}=\{\sqrt{\ell_{i}}\}$
(if necessary, a simple linear transformation can be used to convert the
problem to this form).''

A few remarks are necessary at this moment.

\begin{itemize}
\item Note that in \cite[Th.\ 2]{Gao:07} the meaning of ``$\overline{y}^{\ast
}=(\overline{\varsigma},\overline{\sigma})^{T}$ is a KKT point of
$(\mathcal{P}^{d})$'' is not explained. However, due to the fact that the
constraints of problem $(\mathcal{P}^{d})$ are expressed via $\mathcal{S}_{a}%
$, if $\overline{y}^{\ast}\in\operatorname*{int}\mathcal{S}_{a}$ is a critical
point of $P^{d}$ (that is, $\nabla P^{d}(\overline{y}^{\ast})=0$) then
$\overline{y}^{\ast}$ is a KKT point.

\item It is not clear whether the neighborhood $\mathcal{X}_{o}\times
\mathcal{S}_{o}$ is {}\textquotedblleft a priori\textquotedblright\ prescribed
or the statement should be understood in the sense that there exists such a
neighborhood. In any case the example below shows that \cite[Th.\ 2]{Gao:07}
is false. The proof of this Triality Theorem in \cite{Gao:07} begins with
{}\textquotedblleft In the canonical form of the primal problem (5), replacing
$U(\Lambda(y))$ by the Fenchel-Young equality $(\Lambda(x))^{T}y^{\ast
}-U^{\natural}(y^{\ast})$, the Gao-Strang type \emph{total complementary
function} (see {[}22]) associated with $(\mathcal{P})$ can be obtained as
$\Xi(x,y^{\ast})=\tfrac{1}{2}x^{T}G(\varsigma,\sigma)x-U^{\natural}(y^{\ast
})-x^{T}c-\alpha\varsigma-\ell^{T}\sigma$. $\quad$(21)\textquotedblright.\ For
the proof of the second part of the theorem one says: {}\textquotedblleft On
the other hand, if $\overline{y}^{\ast}\in\mathcal{S}_{a}^{-}$, the matrix
$G(\overline{\varsigma},\overline{\sigma})$ is negative definite. In this
case, the total complementary function $\Xi(x,y^{\ast})$ defined by (21) is a
so-called super-Lagrangian (see {[}12]), i.e., it is locally concave in both
$x\in\mathcal{X}_{o}\subset\mathcal{X}_{a}$ and $y^{\ast}\in\mathcal{S}%
_{o}\subset\mathcal{S}_{a}$. Thus, by the triality theory developed in {[}12],
we have either

\noindent$\displaystyle P(\overline{x})=\min_{x\in\mathcal{X}_{o}}%
P(x)=\min_{x\in\mathcal{X}_{o}}\max_{y^{\ast}\in\mathcal{S}_{o}}\Xi
(x,\lambda)=\min_{y^{\ast}\in\mathcal{S}_{o}}\max_{x\in\mathcal{X}_{o}}%
\Xi(x,\lambda)=\min_{y^{\ast}\in\mathcal{S}_{o}}P^{d}(y^{\ast}),$

\noindent or

\noindent$\displaystyle
P(\overline{x})=\max_{x\in\mathcal{X}_{o}}P(x)=\max_{x\in\mathcal{X}_{o}}%
\max_{y^{\ast}\in\mathcal{S}_{o}}\Xi(x,\lambda)=\max_{y^{\ast}\in
\mathcal{S}_{o}}\max_{x\in\mathcal{X}_{o}}\Xi(x,\lambda)=\max_{y^{\ast}%
\in\mathcal{S}_{o}}P^{d}(y^{\ast}).$

\noindent This proves the statements (19) and (20).\textquotedblright

The references [22] and [12] mentioned above are our references
\cite{Gao/Strang:89} and \cite{Gao-book}, respectively. \newline Therefore the
second part of the conclusion for \cite[Th.\ 2]{Gao:07} does not follow from a
specific results with assumptions that can be verified but from {}%
\textquotedblleft the triality theory\textquotedblright.
\end{itemize}

\begin{example}
\label{ex-th2-g07}Let $n=2$, $A=-4I_{2}$, $B=I_{2}$, $c=(-2,-2)^{T},$
$\alpha=3$, $\ell=(4,4)^{T}$. We have that
\[
P(s,t)=-2s^{2}-2t^{2}+2s+2t+\tfrac{1}{2}\left(  \tfrac{1}{2}s^{2}+\tfrac{1}%
{2}t^{2}-3\right)  ^{2},
\]
and the restrictions are $s^{2}\leq4$, $t^{2}\leq4$, that is $\mathcal{X}%
_{a}=[-2,2]^{2}$. Also,%
\[
P^{d}((y,\sigma,\tau)^{T})=-\frac{2}{y-4+2\sigma}-\frac{2}{y-4+2\tau}%
-\tfrac{1}{2}y^{2}-3y-4\sigma-4\tau.
\]
Then $\overline{y}^{\ast}=(1,1,1)^{T}\in\operatorname*{int}\mathcal{S}_{a}$
and $\overline{y}^{\ast}\in\mathcal{S}_{a}^{-}$ since $G((1,1,1)^{T})=-I_{2}$,
$\overline{y}^{\ast}$ is a KKT point of $P^{d}$ because $\nabla P^{d}%
((1,1,1)^{T})=0$ and $\overline{y}^{\ast}\in\operatorname*{int}\mathcal{S}%
_{a}$, and $\overline{x}=[G((1,1,1)^{T})]^{-1}c=-c=(2,2)^{T}\in\mathcal{X}%
_{a}$. Note that $P(\overline{x})=P^{d}((1,1,1)^{T})=-15/2$. On one hand, for
$\gamma\in(0,1)$ we have that $(2-\gamma,2-\gamma)^{T}\in\mathcal{X}_{a}$ and
\[
P((2-\gamma,2-\gamma)^{T})=-\frac{15}{2}+\frac{1}{2}\gamma^{4}-4\gamma
^{3}+5\gamma^{2}+8\gamma>P(\overline{x}),
\]
which shows that $\overline{x}$ is not a maximum point of $P$ on any
neighborhood of $\overline{x}\in\mathcal{X}_{a}$. Hence relation (20) in the
above theorem does not hold.

On the other hand, for $\gamma\in(0,1)$ we have that
\[
P^{d}((1-16\gamma,1+7\gamma,1+7\gamma)^{T})=-\frac{15}{2}-16\frac{\gamma^{2}%
}{2\gamma+1}\left(  16\gamma+7\right)  <P^{d}((1,1,1)^{T}),
\]
which shows that $\overline{y}^{*}\in\operatorname*{int}\mathcal{S}_{a}$ is
not a local minimum point of $P^{d}$. Hence relation (19) in the above theorem
does not hold, too. Therefore, \cite[Th. 2]{Gao:07} is false.
\end{example}

\section{On a theorem in \cite{FGSW}}

Reference \cite{FGSW} begins with: {}``In this paper, we consider a simple
$0$-$1$ quadratic programming problem in the following form:

$(\mathcal{P})$ : $\min/\max\{P(x)=\tfrac{1}{2}x^{T}Qx-f^{T}x\mid
x\in\mathcal{X}_{a}\}$,$\quad$ (1)

\noindent where $x$ and $f$ are real $n$-vectors, $Q\in\mathbb{R}^{n\times n}$
is a symmetrical matrix of order $n$ and

$\mathcal{X}_{a}=\{x\in\mathbb{R}^{n}\mid0\leq x_{i}\leq1,\ i=1,2,\ldots
,n\}\cap\mathcal{I}^{n}$. $\quad$(2)

\noindent with $\mathcal{I}^{n}=\{x\in\mathbb{R}^{n}\mid x_{i}$ is an integer,
$i=1,2,\ldots,n\}$\textquotedblright, continued with {}\textquotedblleft By
the definition of $\Lambda(x)$ and $V^{\natural}(\sigma)$, we have

$\Xi(x,\sigma)=\tfrac{1}{2}x^{T}Q_{d}(\sigma)x-x^{T}(f+\sigma)$,$\quad$ (8)

\noindent where

$Q_{d}(\sigma)=Q+2\operatorname*{Diag}(\sigma)$

\noindent and $\operatorname*{Diag}(\sigma)\in\mathbb{R}^{n\times n}$ ia a
diagonal matrix with $\sigma_{i}$, $i=1,2,\ldots,n$, being its diagonal
elements\textquotedblright\ and

{}``$P^{d}(\sigma)=-\tfrac{1}{2}(f+\sigma)^{T}Q_{d}^{-1}(\sigma)(f+\sigma
)$.$\quad$ (9)''

Moreover, {}``we introduce the following four sets for consideration:

$\mathcal{S}_{\sharp}^{+}=\left\{  \sigma\in\mathbb{R}^{n}\mid\sigma
>0,\ Q_{d}(\sigma)\text{ is positive definite}\right\}  $, $\quad$(22)

$\mathcal{S}_{\sharp}^{-}=\left\{  \sigma\in\mathbb{R}^{n}\mid\sigma
>0,\ Q_{d}(\sigma)\text{ is negative definite}\right\}  $, $\quad$(23)''

\noindent(we omit the other two sets).

{}``Then we have the following result on the global and local optimality conditions:

\bigskip{}

\textbf{Theorem 3.} Let $Q$ be a symmetric matrix and $f\in\mathbb{R}^{n}$.
Assume that $\overline{\sigma}$ is critical point of $P^{d}(\sigma)$ and
$\overline{x}=\left[  Q_{d}(\overline{\sigma})\right]  ^{-1}(f+\overline
{\sigma})$.

(a) If $\overline{\sigma}\in\mathcal{S}_{\sharp}^{+}$, then $\overline{x}$ is
a global minimizer of $P(x)$ over $\mathcal{X}_{a}$ and $\overline{\sigma}$ is
a global maximizer of $P^{d}(\sigma)$ over $\mathcal{S}_{\sharp}^{+}$ with

$P(\overline{x})=\min\limits_{x\in\mathcal{X}_{a}}P(x)=\max\limits_{\sigma
\in\mathcal{S}_{\sharp}^{+}}P^{d}(\sigma)=P^{d}(\overline{\sigma})$.$\quad$ (26)

(b) If $\overline{\sigma}\in\mathcal{S}_{\sharp}^{-}$, then $\overline{x}$ is
a local minimizer of $P(x)$ over $\mathcal{X}_{a}$ if and only if
$\overline{\sigma}$ is a local minimizer of $P^{d}(\sigma)$ over
$\mathcal{S}_{\sharp}^{-}$, i.e., in a neighborhood $\mathcal{X}_{o}%
\times\mathcal{S}_{o}\subset\mathcal{X}_{a}\times\mathcal{S}_{\sharp}^{-}$ of
$(\overline{x},\overline{\sigma})$,

$P(\overline{x})=\min\limits_{x\in\mathcal{X}_{o}}P(x)=\min\limits_{\sigma
\in\mathcal{S}_{o}}P^{d}(\sigma)=P^{d}(\overline{\sigma})$.$\quad$ (27)''

\bigskip{}

Note that because $\mathcal{X}_{a}$ is a discrete set any $x\in\mathcal{X}%
_{a}$ is a local minimum point for $P$ on $\mathcal{X}_{a}$, as well as a
local maximum point of $P$. In fact the following stronger statement is true.

\begin{theorem}
Let $Q$ be a symmetric matrix and $f\in\mathbb{R}^{n}$. Assume that
$\overline{\sigma}$ is critical point of $P^{d}$ such that $\det
Q_{d}(\overline{\sigma})\neq0$, and $\overline{x}:=\left[  Q_{d}%
(\overline{\sigma})\right]  ^{-1}(f+\overline{\sigma})$. Then $\overline{x}%
\in\mathcal{X}_{a}$ and $P(\overline{x})=\Xi(\overline{x},\overline{\sigma
})=P^{d}(\overline{\sigma}).$

\emph{(a)} If $\overline{\sigma}\in\mathcal{S}_{\sharp}^{+}$, then
$\overline{\sigma}$ is a global maximizer of $P^{d}$ over $\mathcal{S}%
_{\sharp}^{+}$ and $\overline{x}$ is a global minimizer of $P$ over
$\mathcal{X}:=[0,1]^{n};$ in particular, $\overline{x}$ is a global minimizer
of $P$ over $\mathcal{X}_{a}=\{0,1\}^{n}$.

\emph{(b)} If $\overline{\sigma}\in\mathcal{S}_{\sharp}^{-}$, then
$\overline{x}$ is a local minimizer of $P$ over $\mathcal{X}$ and
$\overline{\sigma}$ is a global minimizer of $P^{d}$ over $\mathcal{S}%
_{\sharp}^{-}$.
\end{theorem}

Note that the first part of the above theorem practically covers Theorems 1
and 2 in \cite{FGSW}.

\strut

Proof. It is obvious that $\Xi(x,\cdot)$ is affine (hence concave and convex)
for every $x\in\mathbb{R}^{n}$, $\Xi(\cdot,\sigma)$ is convex for $\sigma
\in\mathcal{S}_{\sharp}^{+}$, and $\Xi(\cdot,\sigma)$ is concave for
$\sigma\in\mathcal{S}_{\sharp}^{-}$. Note that
\begin{equation}
\nabla_{x}\Xi(x,\sigma)=Q_{d}(\sigma)x-(f+\sigma),\quad\nabla_{\sigma}%
\Xi(x,\sigma)(v)=x^{T}\operatorname*{Diag}(v)x-x^{T}v\ \forall v\in
\mathbb{R}^{n}; \label{r3}%
\end{equation}
it follows that $\nabla_{\sigma}\Xi(x,\sigma)=0$ if and only if $x_{i}%
^{2}-x_{i}=0$ for every $i\in\overline{1,n}$, that is, $x\in\mathcal{X}_{a}$.
Furthermore, due to the fact that a critical point of a convex function is a
global minimum point, we have
\begin{equation}
P^{d}(\sigma)=\Xi([Q_{d}(\sigma)]^{-1}(f+\sigma),\sigma)=\left\{
\begin{array}
[c]{ccc}%
\min_{x\in\mathbb{R}^{n}}\Xi(x,\sigma) & \text{if} & \sigma\in\mathcal{S}%
_{\sharp}^{+},\\
\max_{x\in\mathbb{R}^{n}}\Xi(x,\sigma) & \text{if} & \sigma\in\mathcal{S}%
_{\sharp}^{-}.
\end{array}
\right.  \label{r5}%
\end{equation}
Recall the fact that the operator $\varphi:\{U\in\mathfrak{M}_{n}\mid
U\mathrm{\ invertible}\}\rightarrow\mathfrak{M}_{n}$ defined by $\varphi
(U)=U^{-1}$ is Fréchet differentiable and $d\varphi(U)(S)=-U^{-1}SU^{-1}$ for
$U,S\in\mathbb{R}^{n\times n}$ with $U$ invertible, where $\mathfrak{M}_{n}$
is the (normed) linear space of $n\times n$ real matrices. Also, we have
$dQ_{d}(\sigma)(v)=2\operatorname*{Diag}(v)$ and so, on $\mathcal{S}%
^{a}=\{\sigma\in\mathbb{R}^{n}\mid\det Q_{d}(\sigma)\neq0\}$, $d\left[
Q_{d}(\sigma)\right]  ^{-1}(v)=-2\left[  Q_{d}(\sigma)\right]  ^{-1}%
\operatorname*{Diag}(v)\left[  Q_{d}(\sigma)\right]  ^{-1}$ and
\begin{align}
dP^{d}(\sigma)(v) =  &  -v^{T}\left[  Q_{d}(\sigma)\right]  ^{-1}%
(f+\sigma)+(f+\sigma)^{T}\left[  Q_{d}(\sigma)\right]  ^{-1}%
\operatorname*{Diag}(v)\left[  Q_{d}(\sigma)\right]  ^{-1}(f+\sigma
),\label{r1}\\
d^{2}P^{d}(\sigma)(v,v) =  &  -v^{T}\left[  Q_{d}(\mathcal{\sigma})\right]
^{-1}v+4v^{T}\left[  Q_{d}(\mathcal{\sigma})\right]  ^{-1}\operatorname*{Diag}%
(v)\left[  Q_{d}(\mathcal{\sigma})\right]  ^{-1}(f+\sigma)\nonumber\\
&  -4(f+\sigma)^{T}\left[  Q_{d}(\mathcal{\sigma})\right]  ^{-1}%
\operatorname*{Diag}(v)\left[  Q_{d}(\mathcal{\sigma})\right]  ^{-1}%
\operatorname*{Diag}(v)\left[  Q_{d}(\mathcal{\sigma})\right]  ^{-1}(f+\sigma)
\label{r2}%
\end{align}
for all $v\in\mathbb{R}^{n}$.

Since $\overline{\sigma}\in\mathcal{S}^{a}$ is a critical point of $P^{d}$ we
have that $dP^{d}(\overline{\sigma})=0$. Taking into account (\ref{r3}), we
obtain from (\ref{r1}), using a direct computation, that $\nabla_{\sigma}%
\Xi(\overline{x},\overline{\sigma})=0$, and so $\overline{x}\in\mathcal{X}%
_{a}\subset\mathcal{X}.$

Moreover, since $x_{i}^{2}=x_{i}$
\begin{align*}
P(\overline{x})  &  =\tfrac{1}{2}\overline{x}^{T}Q\overline{x}-f^{T}%
\overline{x}=\tfrac{1}{2}\overline{x}^{T}Q_{d}(\overline{\sigma})\overline
{x}-\overline{x}^{T}f-\overline{x}^{T}\operatorname*{Diag}(\overline{\sigma
})\overline{x}\\
&  =\tfrac{1}{2}\overline{x}^{T}Q_{d}(\overline{\sigma})\overline{x}%
-\overline{x}^{T}f-\overline{x}^{T}\overline{\sigma}=\Xi(\overline
{x},\overline{\sigma})\\
&  =\tfrac{1}{2}(f+\overline{\sigma})^{T}[Q_{d}(\overline{\sigma}%
)]^{-1}(f+\overline{\sigma})-\overline{x}^{T}(f+\overline{\sigma}%
)=P^{d}(\overline{\sigma}).
\end{align*}

It is clear that $\mathcal{S}_{\sharp}^{+}$ and $\mathcal{S}_{\sharp}^{-}$ are
open convex sets because $\sigma\rightarrow Q_{d}(\sigma)$ is affine.

If $A:=\left[  Q_{d}(\sigma)\right]  ^{-1}$ is positive definite, setting
$w:=\operatorname*{Diag}(v)\left[  Q_{d}(\mathcal{\sigma})\right]
^{-1}(f+\sigma)$ we have for every $v\in\mathbb{R}^{n}$ that
\[
d^{2}P^{d}(\sigma)(v,v)=-v^{T}Av+4v^{T}Aw-4w^{T}Aw=-(v-2w)^{T}A(v-2w)\leq0,
\]
i.e. $d^{2}P^{d}(\sigma)$ is seminegatively definite. Hence $P^{d}$ is concave
on $\mathcal{S}_{\sharp}^{+}$. Similarly, $P^{d}$ is convex on $\mathcal{S}%
_{\sharp}^{-}.$

(a) Let $\overline{\sigma}\in\mathcal{S}_{\sharp}^{+}$. Since $\Xi
(x,\sigma)=P(x)+\sum_{i=1}^{n}\sigma_{i}(x_{i}^{2}-x_{i})\le P(x)$, for every
$\sigma\ge0$, $x\in[0,1]^{n}$ and taking (\ref{r5}) into account we get%
\begin{align*}
P^{d}(\overline{\sigma})  &  \le\sup_{\sigma\in S_{\sharp}^{+}}P^{d}%
(\sigma)=\sup_{\sigma\in S_{\sharp}^{+}}\min_{x\in\mathbb{R}^{n}}\Xi
(x,\sigma)\le\sup_{\sigma\ge0}\inf_{x\in[0,1]^{n}}\Xi(x,\sigma)\\
&  \le\inf_{x\in[0,1]^{n}}\sup_{\sigma\ge0}\Xi(x,\sigma)\le\inf_{x\in
[0,1]^{n}}P(x)\le P(\overline{x}).
\end{align*}

Therefore $P^{d}(\overline{\sigma})=\max_{\sigma\in\mathcal{S}_{\sharp}^{+}%
}P^{d}(\sigma)$ and $P(\overline{x})=\min_{x\in[0,1]^{n}}P(x)$, since
$P(\overline{x})=P^{d}(\overline{\sigma})$.

(b) Take now $\overline{\sigma}\in\mathcal{S}_{\sharp}^{-}$. Since $P^{d}$ is
convex on $\mathcal{S}_{\sharp}^{-}$ and $\overline{\sigma}$ is a critical
point of $P^{d}$, clearly $\overline{\sigma}$ is a global minimizer of $P^{d}$
on $\mathcal{S}_{\sharp}^{-}.$

Consider $x\in\mathcal{X}$. Since $Q\overline{x}-f=\overline{\sigma
}-2\operatorname*{Diag}(\overline{\sigma})\overline{x}$, we get
\begin{align*}
P(x)  &  =\tfrac{1}{2}x^{T}Qx-f^{T}x=\tfrac{1}{2}(x-\overline{x}%
)^{T}Q(x-\overline{x})+(x-\overline{x})^{T}Q\overline{x}+\tfrac{1}{2}%
\overline{x}^{T}Q\overline{x}-f^{T}x\\
&  =P(\overline{x})+\tfrac{1}{2}(x-\overline{x})^{T}Q(x-\overline
{x})+(x-\overline{x})^{T}Q\overline{x}-(x-\overline{x})^{T}f\\
&  =P(\overline{x})+\tfrac{1}{2}(x-\overline{x})^{T}Q(x-\overline{x}%
)+\sum_{i=1}^{n}\overline{\sigma}_{i}(1-2\overline{x}_{i})(x_{i}-\overline
{x}_{i})=P(\overline{x})+\sum_{i=1}^{n}\mu_{i}(x_{i}-\overline{x}_{i}),
\end{align*}
where $Q=\{q_{ij}\}$ and $\mu_{i}:=\overline{\sigma}_{i}(1-2\overline{x}%
_{i})+\tfrac{1}{2}\sum_{j=1}^{n}q_{ij}(x_{j}-\overline{x}_{j}).$

Let $\varepsilon>0$ be such that $\min_{k\in\overline{1,n}}\overline{\sigma
}_{k}\geq\frac{n}{2}\varepsilon\max_{i,j\in\overline{1,n}}|q_{ij}|$. Take
$U=\{x\in\mathbb{R}^{n}\mid|x_{i}-\overline{x}_{i}|\leq\varepsilon\ \forall
i\in\overline{1,n}\}$. Then $\tfrac{1}{2}\sum_{j=1}^{n}q_{ij}(x_{j}%
-\overline{x}_{j})(x_{i}-\overline{x}_{i})\geq-\frac{n}{2}\varepsilon
|x_{i}-\overline{x}_{i}|\max_{i,j\in\overline{1,n}}|q_{ij}|$ for every
$i\in\overline{1,n}$ and $x\in U$, while the inequality $\overline{\sigma}%
_{i}(1-2\overline{x}_{i})(x_{i}-\overline{x}_{i})\geq|x_{i}-\overline{x}%
_{i}|\min_{k\in\overline{1,n}}\overline{\sigma}_{k}$ for every $i\in
\overline{1,n}$ and $x\in U\cap\mathcal{X}$ is easily checked, since
$\overline{x}_{i}\in\{0,1\}$. This shows that $\mu_{i}(x_{i}-\overline{x}%
_{i})\geq0$ for every $i\in\overline{1,n}$, whence $P(x)\geq P(\overline{x})$
for every $x\in U\cap\mathcal{X}$; therefore, $\overline{x}$ is a local
minimizer of $P$ on $\mathcal{X}.$

\end{document}